\documentclass{ifacconf}

\usepackage{graphicx}      
\usepackage{natbib}        
\usepackage{amsmath,amssymb,mathtools}
\usepackage{tikz}
\usetikzlibrary{shapes,arrows,circuits,decorations.pathmorphing,patterns,positioning,decorations.markings,calc}
\tikzstyle{block} = [draw, fill=blue!20, rectangle, minimum height=3em, minimum width=6em]
\tikzstyle{sum} = [draw, fill=blue!20, circle, node distance=1cm]
\tikzstyle{input} = [coordinate]
\tikzstyle{output} = [coordinate]
\tikzstyle{pinstyle} = [pin edge={to-,thin,black}]

\newtheorem{assumption}{Assumption}
\begin{document}
	\begin{frontmatter}
		
		\title{Strictly negative imaginary state feedback control for relative degree two systems. \thanksref{footnoteinfo}} 
		
		\thanks[footnoteinfo]{This work was supported by the Australian Research Council under grants DP160101121 and DP190102158.}
		
		\author[First]{James Dannatt} 
		\author[Second]{Ian R. Petersen} 
		
		\address[First]{Research School of Electrical, Energy and Materials Engineering, Australian National University, Canberra ACT 2600, Australia (e-mail: jameshd1987@gmail.com).}
		\address[Second]{Research School of Electrical, Energy and Materials Engineering, Australian National University, Canberra ACT 2600, Australia (e-mail: i.r.petersen@gmail.com)}
		
		\begin{abstract}                
			In this paper, we present a strictly negative imaginary state feedback control methodology for relative degree two negative imaginary systems such as flexible structures with collocated sensors and actuators. We show that by augmenting a nominal model with a destabilising PID controller, we can apply state feedback control to the system, resulting in a strictly negative imaginary nominal closed loop system with a prescribed degree of stability. Moreover, the overall system is guaranteed to be internally stable in the face of unmodelled spillover dynamics.
		\end{abstract}
		
		
	\end{frontmatter}
	
	\section{Introduction}
	
	The theory of negative imaginary systems is broadly applicable to problems of robust vibration control for flexible structures \cite{Lanzon2008,petersen2016negative}. These systems are commonly modelled using a sum of second order transfer functions that correspond to the resonant modes of the system \cite{preumont1997vibration,Lanzon2008,reyes2019controller}. The sum can then be split to reduce the order of a system model and the unmodelled dynamics are treated as model uncertainty. When designing a control system for such a reduced order system model, the unmodelled spillover dynamics can degrade control system performance or lead to instability if the controller is not designed to be robust against uncertainty. Negative imaginary systems theory provides a robust way of designing controllers for such systems \cite{Petersen2010,Petersen2015}. 
	
	A key result in NI systems theory states that the positive feedback interconnection of an NI system with a strictly negative imaginary (SNI) system is internally stable when certain gain conditions are satisfied \cite{Lanzon2008,lanzon2017feedback}. A significant advantage of this result is that it guarantees the internal stability of the feedback interconnection through phase stabilisation, rather than the gain stabilisation found in the small-gain theorem. In a small-gain approach, the frequency response of the system would need to satisfy a gain condition across an entire frequency range in order to ensure stability \cite{Petersen2010}. However, negative imaginary systems theory allows us to guarantee stability with only a gain condition applied at DC \cite{lanzon2017feedback}.
	
	Modelling flexible structures with collocated sensors and actuators with a reduced order model is one way to leverage the robust stability properties of NI systems by recognising that the unmodelled dynamics in this approach are guaranteed to have the negative imaginary property \cite{Petersen2010,lanzon2017feedback}. Due to this property, the robust stability of NI systems under a positive feedback interconnection has motivated controller synthesis results with the aim of creating a closed-loop system with the NI or SNI property \cite{Petersen2010,reyes2019controller,Mabrok2015,dannatt2020strictly,REN2021109157,Song2012,10.1016/j.automatica.2022.110235,10.1016/j.automatica.2019.108735}. This closed-loop NI property then guarantees robust stability of the closed-loop to the unmodelled dynamics. Unfortunately, by using a modal sum method of modelling flexible structures with collocated sensors and actuators makes, it can be challenging to apply many of these state feedback results found in the NI literature. The second order transfer functions used in the sum are relative degree two, and due to constraints on the parameters that result in the NI property, these transfer functions will always result in a state space realisation with orthogonal input and output matrices. This problem has been acknowledged in the literature, with authors offering several alternative control schemes; for a comprehensive summary of this problem see Section 7.2 in \cite{reyes2019controller}.
	
	This paper is concerned with new methods for applying static state feedback to a relative degree two system with orthogonal input and output matrices, in order to achieve a closed-loop system which is SNI and the overall system is internally stable. In order overcome the aforementioned orthogonality problem, an augmented control scheme is proposed that modifies the relative degree of the system while giving a designer flexibility in choosing the prescribed degree of stability the closed loop system can achieve. This is important, as this augmented control scheme allows the application of many NI state feedback results from the literature to be applied to systems that previously did not satisfy the orthogonality assumptions these results commonly depend on.
	
\section{Preliminary Results}
	
	In this section, we present some definitions and preliminary results which will be required for our main result. First, consider the linear time-invariant (LTI) system
	\begin{align}
		\dot{{x}}(t) = A{x}(t) + B{u}(t),  \nonumber\\
		{y}(t) = C{x}(t) + D{u}(t). \label{system: LTI system}
	\end{align}
	
	The following definitions relate to the NI, strictly NI and lossless NI properties of the transfer function matrix $G(s) = C(sI-A)^{-1}B + D$, corresponding to the system (\ref{system: LTI system}).
	
	\begin{defn}\cite{Mabrok2015} \label{def:NI} A square transfer function matrix $G(s)$ is NI if the following conditions are satisfied: 
		\begin{enumerate}
			\item G(s) has no pole in $Re[s]>0$.
			
			\item For all $\omega \geq 0$ such that $jw$ is not a pole of $G(s)$, $j(G(j\omega) - G(j\omega)^*) \geq 0$.
			
			\item If $s=j\omega_0$, $ \omega_0 > 0$ is a pole of $G(s)$ then it is a simple pole. Furthermore, if $s=j\omega_0$, $ \omega_0 > 0$ is a pole of $G(s)$, then the residual matrix $K = \lim_{s \to j\omega_0} (s-j\omega_0)jG(s)$ is positive semidefinite Hermitian.
			
			\item If $s=0$ is a pole of $G(s)$, then it is either a simple pole or a double pole. If it is a double pole, then, $\lim_{s \to 0} s^2G(s) \geq 0$.
		\end{enumerate}
	\end{defn}
	
	\begin{defn} \label{def:SNI} \cite{Mabrok2015}
		{A square transfer function matrix $G(s)$ is SNI if the following conditions are satisfied:}
		\begin{enumerate}
			\item G(s) has no poles in $Re[s] \geq 0$.
			\item For all $\omega > 0$ such that $j(G(j\omega) - G(j\omega)^*) > 0$.
		\end{enumerate}
	\end{defn}

	Also, an LTI system (\ref{system: LTI system}) is said to be NI or SNI if the corresponding transfer function matrix $G(s) = C(sI-A)^{-1}B + D$ is NI or SNI.

\begin{defn}\cite{LIU201647} \label{def: LNI definition}
	A square, real, rational, transfer function matrix $G(s)$ is lossless negative imaginary (LNI) if $G(s)$ is negative imaginary and $j[G(j\omega) - G^{\sim}(j\omega)]=0$ for all $\omega>0$ except values of $\omega$ where $j\omega$ is a pole of $G(s)$.
\end{defn}
	
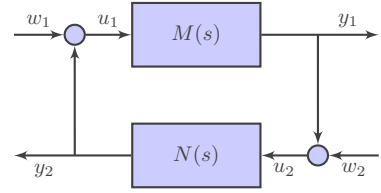
\begin{figure}[!ht]
		\centering
		\begin{tikzpicture}[auto, node distance=2cm,>=latex', black!75,thick, scale=0.8, every node/.style={transform shape}]
			\node [input, name=w1] (w1){};
			\node [sum, right of=w1] (sum1) {};
			\node [block, right of=sum1] (M){$M(s)$};
			\node [coordinate, right of=M] (y1) {};
			\node [output, right of=y1, node distance = 1cm] (y11)  {};
			
			\node [sum, below of=y1, node distance = 2cm] (sum2) {};
			\node [input, right of=sum2, node distance = 1cm] (w2) {};
			\node [block, left of=sum2] (N) {$N(s)$};
			\node [coordinate, left of=N]  (y2) {};
			\node [output, left of=y2, node distance = 1cm] (y22) {};

			\draw [->] (w1) -- node {$w_1$} node [distance=0 and -2cm] {$$} (sum1);
			\draw [->] (sum1) -- node {$u_1$} (M);
			\draw [->] (M) -- node {}(y1) -> node {$y_1$}(y11);
			\draw [->] (y1) -- node {}(sum2) ;
			\draw [->] (w2) -> node {$w_2$}(sum2);
			\draw [->] (sum2) -> node {$u_2$}(N);
			\draw [->] (N) -- node {}(y2) -> node {$y_2$}(y22);
			\draw [->] (y2) -> node {}(sum1) ;
			\node [coordinate, below left of=sum1, node distance = 1cm] {$-$};
			
		\end{tikzpicture}
		\caption{A positive feedback interconnection. The transfer function $M(s)$ and $N(s)$ are interconnected by positive feedback. \cite{Petersen2010}}
		\label{prelim: fig: positive feedback interconnection}
	\end{figure}

\begin{lem} \cite{lanzon2017feedback} \label{lanzon 2017 stability result} \\
		Let $M(s)$ be an NI transfer function matrix without poles at the origin and $N(S)$ be an SNI transfer function matrix. Then the positive-feedback interconnection of $M(s)$ with $N(s)$ as shown in Figure~\ref{prelim: fig: positive feedback interconnection}, is internally stable if and only if the following conditions are satisfied:
		\begin{enumerate}
			\item $I - M(\infty)N(\infty)$ is nonsingular,
			\item $\bar{\lambda}[(I - M(\infty)N(\infty))^{-1}(M(\infty)N(0)-I)]<0$, 
			\item $\bar{\lambda}[(I - N(0)M(\infty))^{-1}(M(0)N(0)-I)]<0$,
		\end{enumerate}
		where $\bar{\lambda}[H]$ denotes the largest eigenvalue of a square, real-valued complex matrix $H$.
\end{lem}
	
\begin{figure}[!ht]
	\centering
	\begin{tikzpicture}[auto, node distance=2cm,>=latex']
		\node [block, name=uncertainty] {$\Delta(s)$};
		\node [block, below of=uncertainty] (plant) {$G(s)$};
		\node [block, below of=plant] (controller) {K};
		\node [coordinate, left of=plant]  (input) {};
		\node [coordinate, right of=plant] (output) {};
		
		\draw [->] (uncertainty) -| node[left, yshift=-5ex]{$w$}([yshift= 1ex]input) -- ([yshift= 1ex]plant.west) ;
		\draw [<-] ([yshift= -1ex]plant.west) -- ([yshift= -1ex]input) |- node[left, yshift=5ex]{$u$} (controller);
		\draw [->] ([yshift= 1ex]plant.east) -- ([yshift= 1ex]output) |- node[right, yshift=-5ex]{$z$} (uncertainty);
		\draw [<-] (controller) -| node[right, yshift=5ex]{$x$}([yshift= -1ex]output) -- ([yshift= -1ex]plant.east);
	\end{tikzpicture}
	\caption{A state feedback control system with plant uncertainty $\Delta(s)$.}
	\label{F1}
\end{figure}
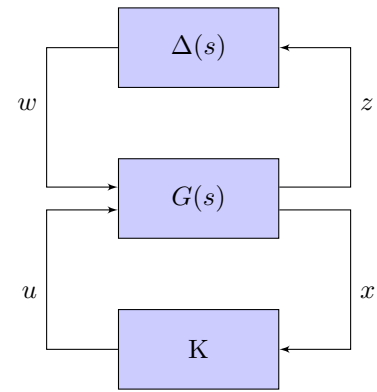

	This paper is concerned with applying state feedback to a relative degree two system with a state space realisation that has orthogonal input and output matrices, in order to achieve a closed-loop system which is SNI with a prescribed degree of stability. Suppose that a given system has enough sensors available that all of a nominal system's states are available for feedback. The system plant can be described by the following linear state space system:
	\begin{align} 
		\dot{{x}} &= A{x} + B_1{w}+B_2{u}, \nonumber\\
		{z} &= C_1{x}. \label{uncertain system}
	\end{align}
	This system corresponds to the nominal plant transfer function matrix $G(s)$ shown in Figure~\ref{F1}. If a state feedback control law $u=Kx$ is applied to this system, the corresponding closed-loop uncertain system has state space representation
	\begin{align}
		\dot{{x}} &= (A + B_2K){x} + B_1{w}, \nonumber\\
		{z}       &= C_1{x}. \label{prelim: math: sfb closed loop system}
	\end{align}
	with corresponding closed-loop transfer function $G_{cl}(s) = C_1(sI - A - B_2K)^{-1}B_1$. Here, we are concerned with finding conditions under which a controller $K$ exists such that the closed-loop transfer function (\ref{prelim: math: sfb closed loop system}) is SNI and internally stable. If a controller exists such that $G_{cl}(s) = C_1(sI - A - B_2K)^{-1}B_1$ is NI (SNI) and the uncertainty $\Delta(s)$ is SNI (NI), then the associated linear fractional interconnection from $w$ to $z$ shown in Figure~\ref{F1} is guaranteed to be internally stable when the conditions of Lemma~\ref{lanzon 2017 stability result} are satisfied. 
	
	The results of this paper build on the state feedback results of the papers \cite{Mabrok2012a,Mabrok2015,dannatt2020strictly}. These state feedback results assume that the system (\ref{uncertain system}) satisfies the following assumptions:
	
	\begin{assumption}
		\label{A1}
		The matrix $C_1B_2$ is non-singular.
	\end{assumption}
	\begin{assumption}
		\label{A2}
		The matrix $R = C_1B_1 + B_1^TC_1^T > 0$.	
	\end{assumption}

	In order to present the state feedback methods proposed in \cite{Mabrok2012a,Mabrok2015,dannatt2020strictly}, we first need to consider the Schur decomposition of the matrix
	\begin{align} 
		A_r &= Q(A+\epsilon I) = A_q + \epsilon Q, \label{math: A_r equation} 
	\end{align} 
	where
	\begin{align}
		Q &= I-B_2(C_1B_2)^{-1}C_1, \label{math: Schur Q definition}\\
		A_q &= A-B_2(C_1B_2)^{-1}C_1A = QA. \label{math: A_q equation}
	\end{align} 
	and $\epsilon>0$ is a parameter defining the required degree of stability of the closed-loop system. We now construct the following matrices:
	\begin{subequations} \label{math: Schur decomp equations}
		\begin{align}
			\tilde{A} &= U^TA_rU =
			\begin{bmatrix}
				\begin{tabular}{ l  r }
					$\tilde{A}_{11}$ & $\tilde{A}_{12}$ \\
					0 & $\tilde{A}_{22}$
				\end{tabular}
			\end{bmatrix}, \label{math: AF}\\
			\tilde{B} &= U^TB_1 =
			\begin{bmatrix}
				\begin{tabular}{c}
					$\tilde{B}_{11}$ \\
					$\tilde{B}_{22}$
				\end{tabular}
			\end{bmatrix}, \label{math: B1}\\
			\tilde{C} &= U^T\big(B_2(C_1B_2)^{-1} - B_1R^{-1}\big) =
			\begin{bmatrix}
				\begin{tabular}{c}
					$\tilde{C}_{11}$ \\
					$\tilde{C}_{22}$
				\end{tabular}
			\end{bmatrix}, \label{math: Bf}\\
			\tilde{Z} &= U^TZU = \tilde{B}R^{-1}\tilde{B}^T - \tilde{C}R\tilde{C}^T = \begin{bmatrix}
				\begin{tabular}{ l  r }
					$\tilde{Z}_{11}$ & $\tilde{Z}_{12}$ \\
					$\tilde{Z}_{21}$ & $\tilde{Z}_{22}$
				\end{tabular}
			\end{bmatrix}, \label{math: Zf}
		\end{align}
	\end{subequations}
	where $\tilde{A}_{11}$ has all of its eigenvalues in the closed left half of the complex plane and $\tilde{A}_{22}$ is an anti-stable matrix; i.e., $\sigma(\tilde{A}_{22}) \subset \{s:\operatorname{Re}[s]>0\}$. Here, $U$ is an orthogonal matrix dependant on $\epsilon$ and is obtained through the real Schur transformation; see Section 5.4 of \cite{Bernstein2009}. 
	
	\begin{rem} \label{R1}
		The eigenvalues of $A_q$ correspond to the zeros of the transfer function from $u$ to $z$ in (\ref{uncertain system}), plus a zero at the origin. To see this, set $w\equiv 0$ and note that (\ref{uncertain system}) implies
		\begin{align}
			\dot{z} = C_1Ax+C_1B_2u. \label{pertb ch: math z dot}
		\end{align}
		Since $(C_1B_2)$ is invertible by assumption, we can rearrange (\ref{pertb ch: math z dot}) to write
		\begin{align}
			u = (C_1B_2)^{-1}\dot{z} - (C_1B_2)^{-1}C_1Ax. \label{pertb ch: math u as function of zdot}
		\end{align}
		Now, we can substitute (\ref{pertb ch: math u as function of zdot}) into (\ref{uncertain system}) to give
		\begin{align}
			\dot{x} &= Ax + B_2(C_1B_2)^{-1}\dot{z}- B_2(C_1B_2)^{-1}C_1Ax \nonumber\\
			&= (A- B_2(C_1B_2)^{-1}C_1A)x + B_2(C_1B_2)^{-1}\dot{z} \nonumber\\
			&= A_qx + B_2(C_1B_2)^{-1}\dot{z} \label{pertb ch: math x dot}
		\end{align}
		These equations define the inverse system which maps from $\dot{z}$ to $u$. Considering the Laplace transform of $\dot{z}$, suppose $z(0)=0$ and observe that $\mathcal{L}(\dot{z}) = s\mathcal{L}(z) + z(0)$. Thus, the eigenvalues of $A_q$ will be the zeros of the transfer function from $u$ to $z$ plus a zero at zero.
	\end{rem}

	The theorem that follows is a special case of Theorem 17 in \cite{dannatt2020strictly}.
	
	\begin{thm}\label{theorem: SNI ARE synthesis theorem no A22} Consider the uncertain system (\ref{uncertain system}) with $r=1$ that satisfies Assumptions A1-A2. Suppose the matrix $A_q$, obtained from the Schur decomposition (\ref{math: Schur decomp equations}) has no unstable eigenvalues. For a given $\epsilon > 0$ such that the corresponding matrix $A_r$ also has no unstable eigenvalues, there exists a static state feedback matrix $K$ such that the closed-loop system (\ref{prelim: math: sfb closed loop system}) is SNI with degree of stability $\epsilon$. The required state feedback controller matrix $K$ is given by
		\begin{equation} 
			K = -(C_1B_2)^{-1}(C_1A + \epsilon C_1 ). \label{math: k matrix reduced}
		\end{equation} 
	\end{thm}
\begin{pf}
	Using Lemma 5 in \cite{dannatt2020strictly}, we know that the closed-loop system (\ref{prelim: math: sfb closed loop system}) is NI if there exists a real $P=P^T\geq0$ such that
	\begin{equation}
		P\hat{A} + \hat{A}^TP + PB_1R^{-1}B_1^TP + Q = 0 \label{math: state feedback no a22 ARE}
	\end{equation}
	where
	\begin{align*}
		\hat{A} &= A_{cl} - B_1R^{-1}C_1A_{cl}, &
		R & = C_1B_1 + B_1^TC_1^T, \\
		Q &= A_{cl}^TC_1^TR^{-1}C_1A_{cl},
	\end{align*}
	with $\sigma \bigg(A_{cl}-B_1R^{-1}(CA_{cl}-B^TP)\bigg) \subset \mathbb{C}_{\leq 0}$ and $A_{cl} = A+\epsilon I + B_2K$. Let $P=0$. Then, (\ref{math: state feedback no a22 ARE}) reduces to the condition
	\begin{align}
		A_{cl}^TC_1^TR^{-1}C_1A_{cl} = 0.
	\end{align}
	We will consider the matrix $C_1A_{cl}$. For our choice of $K$, and noting that for a SISO system, $(C_1B_2)$ is a scalar, $C_1A_{cl} = A+\epsilon I + B_2K$ reduces to
	\begin{align}
		C_1A_{cl} &= C_1A+\epsilon C_1 -C_1B_2(C_1B_2)^{-1}(C_1A + \epsilon C_1) \nonumber\\
		&= C_1A+\epsilon C_1 - C_1A - \epsilon C_1 \nonumber\\
		&=0.
	\end{align}
	Therefore, $P=0$ satisfies (\ref{math: state feedback no a22 ARE}). We will now show that $\sigma \bigg(A_{cl}-B_1R^{-1}(CA_{cl}-B^TP)\bigg) \subset \mathbb{C}_{\leq 0}$. 	
	To achieve this, we begin by considering $A_{cl}-B_1R^{-1}(CA_{cl}-B_1^TP)$. When $P=0$, we have
	\begin{align}
		S =& A_{cl}-B_1R^{-1}(CA_{cl}-B_1^TP) \nonumber\\ 
		=& A + \epsilon{I} + B_2K - B_1R^{-1}C_1(A + \epsilon{I}) - B_1R^{-1}C_1B_2K. \nonumber\\
		\shortintertext{Substituting for $K = -(C_1B_2)^{-1}(C_1A + \epsilon C_1 )$ in this expression, we see that}
		S =& A_\epsilon - B_2(C_1B_2)^{-1}C_1A_\epsilon - B_1R^{-1}C_1A_\epsilon + B_1R^{-1}C_1A_\epsilon \nonumber\\
		=& A_\epsilon - B_2(C_1B_2)^{-1}C_1A_\epsilon \nonumber\\
		=& A_r, \label{math: acl = ar}
	\end{align}
	where $A_r$ is defined in (\ref{math: A_r equation}) and $A_\epsilon = A + \epsilon{I}$.
	Thus, (\ref{math: acl = ar}) implies
	\begin{align}
		\sigma(A_{cl}-BR^{-1}(CA_{cl}-B^TP)) =& \sigma(A_r).
	\end{align}
	Since $\sigma(A_r)\subset \mathbb{C}_{\leq 0}$ follows from our assumption that $A_r$ has no unstable eigenvalues, we have $\sigma(A_{cl}-BR^{-1}(CA_{cl}-B^TP)) \subset \mathbb{C}_{\leq 0}$. It then follows from Lemma 5 in \cite{dannatt2020strictly} that the perturbed closed-loop state space realisation is NI with $\sigma(A_{cl}) \subset \mathbb{C}_{\leq 0}$. Theorem 12 in \cite{dannatt2020strictly} then implies that the actual closed-loop system corresponding to the unperturbed system will have all its poles shifted by an amount $\epsilon$ to the left in the complex plane. Consequently, it follows from Theorem 12 in \cite{dannatt2020strictly} that the closed-loop system is SNI with degree of stability $\epsilon$.
\end{pf}

 For the lemma that will follow, consider a matrix $A_q \in \mathbb{R}^{n \times n}$ with distinct eigenvalues $\lambda_1, \cdots , \lambda_{n}$ and $n\geq 2$. Assume that the following conditions hold:
	\begin{enumerate}
		\item For all $i \in \{1,\cdots,n-1\}$, $Re[\lambda_i] < 0$;
		\item For $i = n$, $\lambda_i= 0$.
	\end{enumerate}
	
	Suppose the eigenvalues of $A_q$ are ordered such that
	\begin{align}
		Re[\lambda_1] \leq \cdots \leq Re[\lambda_{n-1}] < 0, \label{math: gamma for no a22 case n x n}
	\end{align}
	and let $\gamma$ be defined as $\gamma = -Re[\lambda_{n-1}]$. In other words, $\gamma$ is the real part of the maximum eigenvalue of $A_q$ less than zero.
	
	\begin{lem} \label{chap how far pert: distinct eigenvalues}
		Consider an uncertain system (\ref{uncertain system}) with dimension $n$, that satisfies the conditions of Theorem 17 in \cite{dannatt2020strictly} when unperturbed. Suppose the following conditions hold:
		\begin{enumerate}
			\item The pair $\{A,B_2\}$ is controllable;
			\item The matrix $A_q$, (\ref{math: A_q equation}), has $n$ distinct eigenvalues;
			\item The matrix $A_q$, has only one unstable eigenvalue;
		\end{enumerate}
		Then, the closed-loop system formed using Theorem 17 in \cite{dannatt2020strictly} will be SNI with degree of stability $\epsilon$, for any $\epsilon$ within the range $0 < \epsilon < \gamma$, where $\gamma = -Re[\lambda_{n-2}]$.
	\end{lem}

	The following theorem outlines the degree of stability that can be attained by the closed-loop system (\ref{prelim: math: sfb closed loop system}) formed using (\ref{math: k matrix reduced}) under state feedback.
	
	\begin{thm}  \label{chap how far pert: distinct eigenvalues no A22}
		Consider the system (\ref{uncertain system}) with dimension $n$, that satisfies the conditions of Theorem \ref{theorem: SNI ARE synthesis theorem no A22} when unperturbed. Suppose, 
		the pair $\{A,B_2\}$ is controllable;
		the matrix $A_q$ defined in  (\ref{math: A_q equation}), has distinct eigenvalues and 
		the matrix $A_q$, has no unstable eigenvalues.
		Then, the closed-loop system formed using Theorem \ref{theorem: SNI ARE synthesis theorem no A22} will be SNI with degree of stability $\epsilon$, for any $\epsilon$ within the range $0 < \epsilon < \gamma$, where $\gamma = -Re[\lambda_{n-1}]$. In addition,  for $\epsilon$ in the range $-Re[\lambda_{n-1}] < \epsilon < -Re[\lambda_{n-2}]$, if the open-loop system satisfies the conditions of Theorem 17 in \cite{dannatt2020strictly} at $\epsilon$, then the system will remain SNI for all $\epsilon$ within the range $-Re[\lambda_{n-2}] < \epsilon < -Re[\lambda_{n-2}]$.
	\end{thm} 
	\begin{pf}
		For $\epsilon$ within the range $0 < \epsilon < \gamma$, where $\gamma = -Re[\lambda_{n-1}]$, the closed-loop system is SNI follows directly from Theorem~\ref{theorem: SNI ARE synthesis theorem no A22}. For $-Re[\lambda_{n-2}] < \epsilon < -Re[\lambda_{n-2}]$, the matrix $A_q$ will have attained an unstable eigenvalue. Consequently, $X$ has dimension one. Therefore, if the system satisfies the conditions of Theorem 17 in \cite{dannatt2020strictly} for $-Re[\lambda_{n-2}] < \epsilon < -Re[\lambda_{n-2}]$, then the system will remain SNI for all $\epsilon$ within the range $0< \epsilon < -Re[\lambda_{n-2}]$ follows directly from Theorem\ref{chap how far pert: distinct eigenvalues}. 
	\end{pf}

\section{Modelling NI systems for a robust control framework} \label{sec 3}

In our introduction, we stated that modelling flexible structures with collocated sensors and actuators using a sum of 2nd order transfer functions that correspond to the system natural modes was convenient as unmodelled dynamics in this approach are guaranteed to have the negative imaginary property. In the following section we outline two methods of modelling NI systems for use within a robust control framework. The first uses the additive uncertainty structure discussed in our introduction. The second, considers a new multiplicative uncertainty structure for lossless NI systems. 

\subsection{Modelling NI systems with Additive uncertainty}
A lightly damped structures with collocated position sensors and force actuators may be modelled using an infinite sum of second-order functions,
	\begin{align}
		G_T(s) = \sum_{p=1}^{M}\frac{\Gamma_p}{s^2 + 2\zeta_p\omega_ps + \omega_p^2}, \label{ch modell: math: assumed modes full model}
	\end{align}
	where $\Gamma_p=\Psi_p \Psi_p^T \geq 0$, $\omega_p$ denotes the natural frequency, $\zeta_p$ is analogous to the damping coefficient, and $\Psi_p$ is an $n \times 1$ vector \cite{preumont1997vibration}. By choosing a small number of modes for a reduced order nominal model, the sum (\ref{ch modell: math: assumed modes full model}) can be split into two sums as 
	\begin{align}
		G_n(s) &= \sum_{p=1}^{m}\frac{\Gamma_p}{s^2 + 2\zeta_p\omega_ps + \omega_p^2}, \label{ch modelling: assumed modes: nom plant} \\
		\Delta(s) &= \sum_{p=m+1}^{M}\frac{\Gamma_p}{s^2 + 2\zeta_p\omega_ps + \omega_p^2}, \label{ch modelling: assumed modes: ucnertainty}
	\end{align}
	where $m$ represents the desired number of modes in the nominal model, $G_n(s)$ represents the truncated 'nominal' model, and $\Delta(s)$ represents the unmodelled system dynamics corresponding to the lightly damped modes not captured in the nominal model. In design, we can now work only with the reduced order nominal model that considers only the system modes in a bandwidth of interest. In this scheme, we consider $\Delta(s)$ as a model uncertainty that is represented by the additive perturbation  
	\begin{align}
		G_T(s) = G_n(s) + \Delta(s). \label{ch model: math additive uncert: G + delta}
	\end{align}
	We can readily verify that the uncertainty $\Delta(s)$ is both stable and has the negative imaginary property according to Definition~\ref{def:NI}. This method of treating unmodelled spillover dynamics in flexible systems with collocated force actuators and position sensors is commonly used, since it is straight forward to verify that the unmodelled dynamics will have the NI property. We will show in Section~\ref{sec 4} that this method does not satisfy the assumptions needed to apply Theorem~\ref{theorem: SNI ARE synthesis theorem no A22} or any NI feedback result that requires the matrices $C_1B_1 > 0$.
	
\subsection{Multiplicative uncertainty}

	Consider now an alternative model for NI systems that can be described using the multiplicative perturbation
	\begin{align}
		G_T(s) &= G_n(s)(1 + \Delta(s)), \label{ch modelling: math: full lossless system}
	\end{align}
	where $G_n(s)$ is our nominal plant and $\Delta(s)\in \mathcal{RH}_{\infty}$ is lossless negative imaginary. In order to create a model (\ref{ch modelling: math: full lossless system}) such that the uncertainty is LNI, suppose a given system can be modelled as an infinite sum of second-order transfer functions as in(\ref{ch modell: math: assumed modes full model}). However, assume also that the model is undamped (conservative), with transfer function matrix
	\begin{align}
		G_T(s) = \sum_{i=1}^{\infty}\frac{\Gamma_i}{s^2+\omega_i^2}, \label{ch modelling: math: multiplicative undamped sum}
	\end{align}
	where $\Gamma_i> 0$ and $\omega_i$ denotes the natural frequency of the '$i^{th}$' mode.
	
	For practical modelling purposes, we choose a finite and often small number of modes to model. This results in a reduced order model of the form
	\begin{align}
		G_n(s) = \sum_{i=1}^{M}\frac{\Gamma_i}{s^2+\omega_i^2},
	\end{align}
	where $M$ is the desired number of modes. Then it follows from (\ref{ch modelling: math: full lossless system}) that
	\begin{align}
		\Delta(s) = \frac{G_T(s)}{G_n(s)} -1 = \frac{G_T(s)-G_n(s)}{G_n(s)} = \frac{\sum_{i=M+1}^{\infty}\frac{\Gamma_i}{s^2+\omega_i^2}}{\sum_{i=1}^{M}\frac{\Gamma_i}{s^2+\omega_i^2}}. \label{ch modelling: math: lossless system uncertainty}
	\end{align}
	
	It is clear from the structure of (\ref{ch modelling: math: lossless system uncertainty}) that $\Delta(j\omega)$ is real for all values of $\omega$ (see equation 2.5 in \cite{preumont1997vibration}) and according to Definition~\ref{def: LNI definition}, is lossless negative imaginary. Similar to using an additive uncertainty structure, this method does not satisfy the assumptions needed to apply Theorem~\ref{theorem: SNI ARE synthesis theorem no A22} or any NI feedback result that requires the matrices $C_1B_1 > 0$.
	
\section{Main result} \label{sec 4}
	
\begin{figure}[!ht]
		\centering
		\begin{tikzpicture}[auto, node distance=2cm,>=latex', black!75,thick, scale=0.6, every node/.style={transform shape}]
			
			\node [block, name=uncertainty] {$\Delta(s)$};
			\node [block, below of=uncertainty] (plant) {$G_n(s)$};
			\node [block, below of=plant] (alpha) {$-1$};
			\node [coordinate, left of=plant]  (nom_input) {};
			\node [block, left of=nom_input]  (integrator) {$\frac{\tilde{K}}{s}$};
			\node [sum, left of=integrator, node distance = 2.2cm]  (input_sum) {+};
			\node [coordinate, left of=input_sum]  (input_tip) {};
			\node [coordinate, above of=input_sum]  (input_u) {};
			\node [coordinate, below of=input_sum]  (input_fb) {};
			\node [sum, right of=plant, node distance = 2.2cm] (sum) {+};
			\node [coordinate, right of=sum, distance = 0.3cm]  (feedback_branch) {};
			\node [coordinate, below of=feedback_branch]  (output_fb) {};	
			\node [coordinate, right of=feedback_branch] (output) {};
			
			\draw [-] (integrator.east) -- node[above, yshift=0ex]{z}(nom_input) ;
			\draw [->] (nom_input) -- node[below, yshift=0ex]{} (plant);
			\draw [->]  (nom_input) |- node[left, yshift=0ex]{} (uncertainty);
			\draw [->] (plant) -- node[above, yshift=0ex]{$\tilde{y}$}(sum) ;
			\draw [<-] (sum) |- node[right, yshift=-5ex]{$w$} (uncertainty);
			\draw [->] (sum) -- node[above, yshift=0ex]{$y$}(output) ;
			\draw [->] (input_sum) -- node[above, yshift=0ex]{} (integrator.west);
			\draw [->] (input_tip) -- node[above, yshift=0ex]{$v$} (input_sum.west);
			\draw [->] (input_u) -- node[left, yshift=0ex]{$u$} (input_sum.north);
			\draw [->] (input_fb) -- (input_sum.south);
			\draw [->] (feedback_branch) -- (output_fb) -> (alpha.east);
			\draw [-] (alpha.west) -- (input_fb);
		\end{tikzpicture}
		\caption{Block diagram of a flexible cantilever model augmented with integral control.}
		\label{chap model: fig: closed-loop model}
	\end{figure}
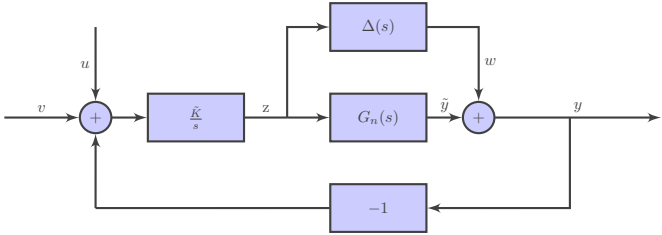

	Suppose we have a system that we wish to control using SNI state feedback. If the system was modelled using either method from Section~\ref{sec 3}, it is straight forward to verify that neither model satisfies the assumptions A1-A2. In either case, due to the structure of the transfer function model, the corresponding state space realisation (\ref{uncertain system}), is such that the output matrix $C_1$ is orthogonal to the input matrices $B_1$ and $B_2$. In this section we propose a method of overcoming this by augmenting the system with a controller that alters the dynamics of the open loop system. To that end, consider the augmented system shown in Figure~\ref{chap model: fig: closed-loop model}. Suppose that an NI system has been modelled using an additive uncertainty structure as in (\ref{ch model: math additive uncert: G + delta}) with only the first mode being considered for control. The nominal plant $G_n(s)$ can be realised with the following state space realisation:
	\begin{align} 
		\dot{\tilde{x}} &= \tilde{A}\tilde{x} + \tilde{B}z, \\
		\tilde{y} &= \tilde{C}\tilde{x},
	\end{align}
	where
	\begin{align} 
		\tilde{A} &= \begin{bmatrix}
			\begin{tabular}{ c c } 	
				$-2\zeta_1\omega_1$ & $-\omega_1^2$ \\ $1$ & $0$
			\end{tabular}
		\end{bmatrix}, & \tilde{B} &= \begin{bmatrix}
			\begin{tabular}{ c } 	
				$1$ \\
				$0$ 
			\end{tabular}
		\end{bmatrix}, \nonumber\\ \tilde{C} &= \begin{bmatrix}
			\begin{tabular}{ c c } 	
				$0$ & $\Gamma_1$ 
			\end{tabular}
		\end{bmatrix}, & \tilde{x} &= \begin{bmatrix}
			\begin{tabular}{ c } 	
				$x_1$ \\  $x_2$
			\end{tabular}
		\end{bmatrix}. \label{ch model: math: nominal system state space realisation}
	\end{align}
	
	In this scheme, the nominal plant is augmented using an integral controller in order to alter the system such that it satisfies the assumptions A1-A2. Using (\ref{ch model: math: nominal system state space realisation}) and Figure~\ref{chap model: fig: closed-loop model}, we can form a realisation for the augmented closed-loop system. Let $z=x_3$ and a state space realisation for the augmented system is given by	
	\begin{align} 
		\dot{x} &= Ax + B_1w + B_2u , \\
		z &= C_1x, \label{ch model: math: augmented system}
	\end{align}
	where
	\begin{align} 
		A &= \begin{bmatrix}
			\begin{tabular}{ c c } 	
				$\tilde{A}$ & $\tilde{B}$ \\ $-\tilde{K}\tilde{C}$ & $0$
			\end{tabular}
		\end{bmatrix}, & B_1 &= \begin{bmatrix}
			\begin{tabular}{ c } 	
				$0$ \\
				$-\tilde{K}$ 
			\end{tabular}
		\end{bmatrix}, &
		B_2 &= \begin{bmatrix}
			\begin{tabular}{ c } 	
				$0$ \\
				$\tilde{K}$ 
			\end{tabular}
		\end{bmatrix}, \nonumber \\ 
		C_1 &= \begin{bmatrix}
			\begin{tabular}{ c c } 	
				$0$ & $1$ 
			\end{tabular}
		\end{bmatrix}, & 
		x &= \begin{bmatrix}
			\begin{tabular}{ c } 	
				$\tilde{x}$ \\  $x_3$
			\end{tabular}
		\end{bmatrix}. & \label{ch model: math: augmented system state space realisation}
	\end{align}
	
	When $\tilde{K} < 0$, this realisation satisfies assumptions A1-A2. Thus, we can freely apply Theorem~\ref{theorem: SNI ARE synthesis theorem no A22} to achieve SNI state feedback using a controller of the form
	\begin{equation} 
		K = -(C_1B_2)^{-1}(C_1A + \epsilon C_1 ) = \begin{bmatrix}
			\begin{tabular}{ c c c} 	
				$0$ & $\Gamma$ & $-\frac{\epsilon}{\tilde{K}}$
			\end{tabular}
		\end{bmatrix}.
	\end{equation} 
	Under the state feedback control law $u=Kx$ the closed-loop system corresponding to the plant (\ref{ch model: math: augmented system}) becomes
	\begin{align} 
		\dot{x} &= \begin{bmatrix}
			\begin{tabular}{ c c c} 	
				$-2\zeta_1\omega_1$ & $-\omega_1^2$ & $1$ \\ 
				$1$ & $0$ & $0$ \\
				$0$ & $0$ & $-\epsilon$
			\end{tabular}
		\end{bmatrix}\begin{bmatrix}
			\begin{tabular}{ c } 	
				$x_1$ \\  $x_2$ \\  $x_3$
			\end{tabular}
		\end{bmatrix} + \begin{bmatrix}
			\begin{tabular}{ c } 	
				$0$ \\  $0$ \\  $-\tilde{K}$
			\end{tabular}
		\end{bmatrix}w , \nonumber\\
		z &= \begin{bmatrix}
			\begin{tabular}{ c c c} 	
				$0$ & $0$ & $1$ 
			\end{tabular}
		\end{bmatrix}\begin{bmatrix}
			\begin{tabular}{ c } 	
				$x_1$ \\  $x_2$ \\  $x_3$
			\end{tabular}
		\end{bmatrix}. \label{ch model: math: augmented system under statefeedback}
	\end{align}
	
	For this system, we can see that the transfer function from $w \to z$ is NI and stable. However, using Theorem~\ref{theorem: SNI ARE synthesis theorem no A22} there is a negligible performance advantage. This is because in this approach, the size of the allowed perturbation is constrained by the location of the open loop poles of the nominal plant. In order to see this, consider the following claim.
	
	\begin{claim} \label{ch modell: example: cant perturb system}
		The degree of stability attainable for the augmented system (\ref{ch model: math: augmented system state space realisation}) using Theorem~\ref{theorem: SNI ARE synthesis theorem no A22} is the degree of stability of the open loop plant, which is $\epsilon = \zeta_1 \omega_1$. Theorem~\ref{chap how far pert: distinct eigenvalues no A22} tells us the largest degree of stability we can achieve using perturbation of an NI system under state feedback. Specifically, the closed-loop system formed using Theorem (\ref{theorem: SNI ARE synthesis theorem no A22}) will be SNI with degree of stability $\epsilon$, for any $\epsilon$ within the range $0 < \epsilon < \gamma$, where $\gamma = -Re[\lambda_{n-1}]$ and $\lambda_{n-1}$ is an eigenvalue of the matrix $A_q=(I-B_2(C_1B_2)^{-1}C_1)A$. We assume the ordered eigenvalues $\lambda_i$ are monotonically increasing. For the augmented system (\ref{ch model: math: augmented system state space realisation}),
		\begin{align}
			A_q = \begin{bmatrix}
				\begin{tabular}{ c c } 	
					$\tilde{A}$ & $\tilde{B}$ \\ $0$ & $0$
				\end{tabular}
			\end{bmatrix}.
		\end{align}
		$A_q$ is upper triangular and therefore the eigenvalues of $A_q$ are the eigenvalues of $\tilde{A}$ and an additional eigenvalue at zero.
		When $\tilde{A}$ is defined as in (\ref{ch model: math: nominal system state space realisation}), the eigenvalues of $\tilde{A}$ are $-\zeta_1\omega_1 \pm \omega_1\sqrt{\zeta_1^2-1}$. Since $A_q$ is upper triangular, we can readily see that 
		\begin{align}
			\sigma(A_q) = \{0, -\zeta_1\omega_1+\omega_1\sqrt{\zeta_1^2-1}, -\zeta_1\omega_1-\omega_1\sqrt{\zeta_1^2-1} \}.
		\end{align}
		Since our system is lightly damped, we can assume that $0\leq \zeta_1 < 1$. Therefore,
		\begin{align}
			\sigma(A_q) = \{0,  -\alpha\pm\beta j\},
		\end{align} 
		where $\alpha = \zeta_1\omega_1\geq 0$ and $\beta = \omega_1\sqrt{1-\zeta_1^2}>0$. According to Theorem~\ref{chap how far pert: distinct eigenvalues no A22}, the closed-loop system formed using Theorem (\ref{theorem: SNI ARE synthesis theorem no A22}) will be SNI with degree of stability $\epsilon$, for any $\epsilon$ within the range $0 < \epsilon < \zeta_1\omega_1$.
	\end{claim}
	
\begin{figure}[!ht]
		\centering
		\begin{tikzpicture}[auto, node distance=2cm,>=latex', black!75,thick, scale=0.6, every node/.style={transform shape}]
			
			\node [block, name=uncertainty] {$\Delta(s)$};
			\node [coordinate, below of=uncertainty] (bunc) {$$};
			\node [block, below of=bunc] (alpha) {$-1$};
			\node [coordinate, left of=bunc]  (nom_input) {};
			\node [block, left of=nom_input]  (integrator) {$G_n(s)$};
			\node [coordinate, left of=integrator]  (plantPIDspacer) {};
			\node [block, left of=plantPIDspacer, node distance = 1cm]  (PID) {$C(s)$};
			\node [sum, left of=PID, node distance = 2.2cm]  (input_sum) {+};
			\node [coordinate, above of=input_sum]  (input_u) {};
			\node [coordinate, below of=input_sum]  (input_fb) {};
			\node [sum, right of=bunc, node distance = 2.2cm] (sum) {+};
			\node [coordinate, right of=sum, node distance = 1cm]  (feedback_branch) {};
			\node [coordinate, below of=feedback_branch]  (output_fb) {};	
			\node [coordinate, right of=feedback_branch, node distance = 1cm] (output) {};
			
			\draw [-] (integrator.east) -- node[above, yshift=0ex]{z}(nom_input) ;
			\draw [->]  (nom_input) |- node[left, yshift=0ex]{} (uncertainty);
			\draw [->] (integrator.east) -- node[above, yshift=0ex]{}(sum) ;
			\draw [<-] (sum) |- node[right, yshift=-5ex]{$w$} (uncertainty);
			\draw [->] (sum) -- node[above, yshift=0ex]{$y$}(output) ;
			\draw [->] (input_sum) -- node[above, yshift=0ex]{} (PID.west);
			\draw [->] (PID.east) -- node[above, yshift=0ex]{} (integrator.west);
			\draw [->] (input_u) -- node[left, yshift=0ex]{$u$} (input_sum.north);
			\draw [->] (input_fb) -- (input_sum.south);
			\draw [->] (feedback_branch) -- (output_fb) -> (alpha.east);
			\draw [-] (alpha.west) -- (input_fb);
		\end{tikzpicture}
		\caption{A flexible cantilever model with lossless NI uncertainty $\Delta(s)$, augmented using a PID controller $C(s)$.}
		\label{chap model: fig: Multiplicative uncertainty under PID with feedback}
	\end{figure}
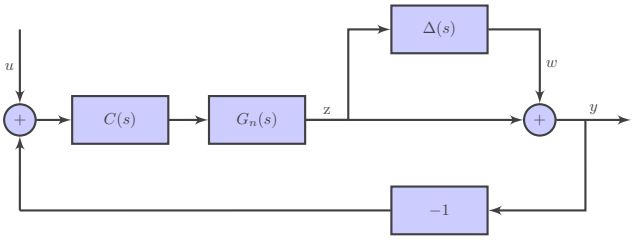

	Augmenting the system model with an integrator was shown to be sufficient for allowing SNI state feedback to be applied to the system. However, the closed loop degree of stability that could be achieved using Theorem~\ref{theorem: SNI ARE synthesis theorem no A22} was constrained by the location of the open loop poles of the nominal plant. In order to overcome this limitation, we now introduce a modified control scheme that augments a nominal system model with a PID controller and allows for the application of SNI state feedback. In addition, we will be able to adjust the attainable degree of stability for the closed-loop SNI system by choice of the PID parameters. 
	Since our approach can be applied to either of the modelling methods from the previous section, we will demonstrate this scheme using a model that has a multiplicative uncertainty structure (\ref{ch modelling: math: multiplicative undamped sum}). Considering only the first mode, we augment the nominal plant using the control scheme shown in Figure~\ref{chap model: fig: Multiplicative uncertainty under PID with feedback}.
	
	For the feedback architecture shown in Figure~\ref{chap model: fig: Multiplicative uncertainty under PID with feedback}, assume that $C(s)$ is a PID controller given by
	\begin{align}
		C(s) &= k_p + \frac{k_i}{s} + k_d{s} = \frac{k_d{s^2} + k_p{s} + k_i}{s}, \label{ch model: math: C(s) }
	\end{align}
	which implies
	\begin{align}
		L(s) = C(s)G_n(s) &= \frac{\Gamma({k_d{s^2} + k_p{s} + k_i})}{s({s^2 + \omega_p^2})} = \frac{n_L(s)}{d_L(s)}. \label{ch model: math: L(s)}
	\end{align}

	\begin{rem}
	In order to see how the choice of controller $C(s)$ affects the degree of stability that can be attained using Theorem~\ref{theorem: SNI ARE synthesis theorem no A22}, we define $T_{uz}(s)$ as the closed-loop transfer function from input $u$ to ouput $z$ and observe that 
	\begin{align}
		T_{uz}(s) = \frac{L(s)}{1+L(s)} = \frac{n_L(s)}{n_L(s)+d_L(s)}.
	\end{align}
	In Remark~\ref{R1}, it was demonstrated that the open loop zeros of $T_{uz}(s)$ correspond to the eigenvalues of (\ref{math: A_q equation}). For this system, the open loop zeros of $T_{uz}(s)$ are the roots of $n_L(s)$. Also, Theorem~\ref{chap how far pert: distinct eigenvalues no A22} outlines that the degree of stability that can be attained by a closed-loop system (\ref{prelim: math: sfb closed loop system}) formed using (\ref{math: k matrix reduced}) under state feedback is dictated by the real part of the eigenvalues of $A_q$. Consequently, the prescribed degree of stability attainable using Theorem~\ref{theorem: SNI ARE synthesis theorem no A22} is determined by the roots of $n_L(s)$.
	\end{rem}
	
	Before we proceed with the control problem, we must check that the augmented system satisfies assumptions A1-A2. A state space realisation for the augmented plant $L(s)=C(s)G_n(s)$ defined in (\ref{ch model: math: L(s)}), in controllable canonical form is given by
	\begin{align}
		\dot{{x}} &= \begin{bmatrix}
			\begin{tabular}{ c c c} 	
				$0$ & $1$ & $0$ \\ 
				$0$ & $0$ & $1$ \\
				$0$ & $-\omega_p^2$ & $0$
			\end{tabular}
		\end{bmatrix}{x} + \begin{bmatrix}
			\begin{tabular}{ c } 	
				$0$ \\
				$0$ \\
				$1$
			\end{tabular}
		\end{bmatrix}e, \\
		z &= \begin{bmatrix}
			\begin{tabular}{ c c c} 	
				$\Gamma{k_i}$ & $\Gamma{k_p}$ & $\Gamma{k_d}$ 
			\end{tabular}
		\end{bmatrix}{x}.
	\end{align}
	
	From Figure~\ref{chap model: fig: Multiplicative uncertainty under PID with feedback} we see that input to the augmented plant, $e$ is equal to $e = u-(w+z)$. Thus, a state space realisation for the closed-loop system shown in Figure~\ref{chap model: fig: Multiplicative uncertainty under PID with feedback} is given by
	\begin{align}
		\dot{x} &= Ax + B_1w + B_2u, \nonumber\\
		z       &= C_1x, 
	\end{align}
	where
	\begin{align}
		A &= \begin{bmatrix}
			\begin{tabular}{ c c c} 	
				$0$ & $1$ & $0$ \\ 
				$0$ & $0$ & $1$ \\
				$-\Gamma{k_i}$ & $-(\Gamma{k_p}+\omega_p^2)$ & $-\Gamma{k_d}$
			\end{tabular}
		\end{bmatrix}, \\ B_2 &= \begin{bmatrix}
			\begin{tabular}{ c } 	
				$0$ \\
				$0$ \\
				$1$
			\end{tabular}
		\end{bmatrix}, \quad B_1 = \begin{bmatrix}
			\begin{tabular}{ c } 	
				$0$ \\
				$0$ \\
				$-1$
			\end{tabular}
		\end{bmatrix}, \\
		C_1 &= \begin{bmatrix}
			\begin{tabular}{ c c c} 	
				$\Gamma{k_i}$ & $\Gamma{k_p}$ & $\Gamma{k_d}$ 
			\end{tabular}
		\end{bmatrix}.
	\end{align}
	Since $\Gamma>0$ by assumption, by choosing ${k_d<0}$ we can readily verify that the SNI state feedback assumptions A1-A2 are satisfied. In addition, by choosing the terms $k_d$, $k_i$ and $k_p$, the degree of stability the closed-loop system can achieve may be arbitrarily chosen. While this seems nonrestrictive, bandwidth requirements and available hardware will limit the realisable PID values that can be used.
	
	\section{Illustrative Example}
	
\begin{figure}[!ht]
		\centering
		\includegraphics[width=3in]{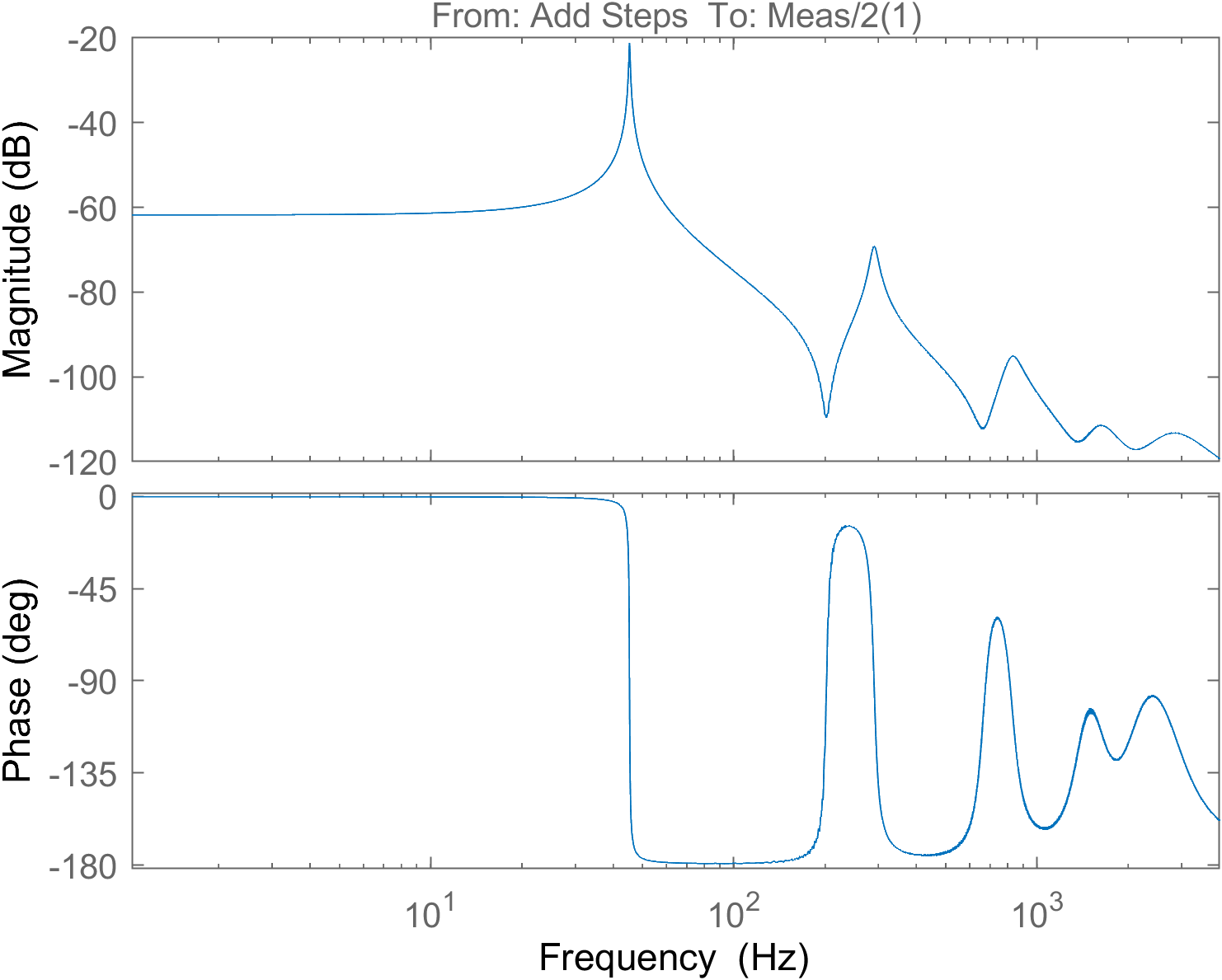}
		\caption{The frequency response of a flexible cantilever, measured from a force input at the tip to the vertical displacement of the tip.}
		\label{ch examp: fig: open loop freq response of beam}
	\end{figure}

	In this section we provide an illustrative example that demonstrates how to apply our main results to a relative degree two system with a corresponding state space realisation that does not satisfy Assumptions 1-2. For this example, we will apply SNI state feedback to a simulated flexible cantilever with the frequency response shown in Figure~\ref{ch examp: fig: open loop freq response of beam}.
	
	The frequency response shown in Figure~\ref{ch examp: fig: open loop freq response of beam} can be fit to a transfer function model of the form
	\begin{align}
		G(s) = \sum_{i=1}^{M}\frac{\Gamma_i}{s^2 + \gamma_i s + \delta_i}, \label{ch example: math: full modal transfer function sum}
	\end{align}
	where $\delta_i=\omega_i^2$. Each $\omega_i>0$ represents a modal frequency, $\gamma_i=2\zeta_i\omega_i$, $\zeta_i$ is loosely the damping term, and $\Gamma_i>0$. To achieve this, we use the linearisation toolbox in Matlab. For this example, only the first 5 modes ($M=5$) are considered, which correspond to the values shown in Table~\ref{Table1}.
	
	\begin{table}[h!]
		\centering
		\begin{tabular}{c | c | c | c}
			i & $\Gamma_i$ & $\gamma_i$ & $\delta_i$ \\ \hline
			1 & $64.06$ & $2.089$ & $8.096\times 10^{4}$ \\
			2 & $65.14$ & $85.83$ & $3.327\times 10^{6}$ \\
			3 & $63.57$ & $695.7$ & $2.697\times 10^{7}$ \\
			4 & $51.76$ & $2624$  & $1.017\times 10^{8}$ \\
			5 & $274$ 	& $8014$  & $3.106\times 10^{8}$
		\end{tabular}
		\caption{Parameters of (\ref{ch example: math: full modal transfer function sum}) for M=5, corresponding to the frequency response in Figure~\ref{ch examp: fig: open loop freq response of beam}.}
		\label{Table1}
	\end{table}
	
	Referring to Table~\ref{Table1}, we can see that the flexible cantilever is lightly damped and the first resonant mode occurs at approximately $45$ Hz. Rather than designing a control scheme for the entire beam, a model reduction is performed and we consider only the first mode. Now, we treat the system is if it is conservative (undamped). Thus, we may approximate the system model by
	\begin{align}
		G_n(s) = \frac{64.06}{s^2+8.096\times10^4}. \label{ch exam: math: lightly damped GN}
	\end{align}
	The remaining system dynamics are treated as a multiplicative uncertainty and were shown in (\ref{ch modelling: math: lossless system uncertainty}) to be lossless negative imaginary. Using the 2nd order approximation (\ref{ch exam: math: lightly damped GN}), we now augment this nominal plant using the PID control scheme proposed in Figure~\ref{chap model: fig: Multiplicative uncertainty under PID with feedback}. 
	
	When choosing the values for the PID controller we refer back to (\ref{ch model: math: C(s) }) and note that the open-loop zeros of the augmented system are given by the roots of
	\begin{align}
		n_L(s) = \Gamma({k_d{s^2} + k_p{s} + k_i}).
	\end{align}
	Remark~\ref{R1} and Theorem~\ref{chap how far pert: distinct eigenvalues no A22}, tell us that the open-loop zeros of the augmented system dictate the prescribed degree of stability we can obtain using SNI state feedback. In particular, Theorem~\ref{chap how far pert: distinct eigenvalues no A22} states that if we choose our PID values such that the open-loop zeroes are all stable, the closed-loop system formed using SNI state feedback will be SNI with degree of stability $\epsilon$, for any $\epsilon$ within the range $0 < \epsilon < \gamma$. Using this information, we will choose the PID values to allow for a degree of stability $\epsilon=10$. This corresponds to a PID controller $C(s)$ given by
	\begin{align}
		C(s) &= \frac{k_d{s^2} + k_p{s} + k_i}{s} = \frac{-0.2{s^2} -7{s} -50}{s}, \label{ch example: PID controller C(s)}
	\end{align}
	where
	\begin{align}
		k_d &= -0.2, & k_p &= -7, & k_i &= -50.
	\end{align}
	These PID values were chosen arbitrarily in order give a prescribed degree of stability $\epsilon=10$. Based on the location of the open-loop system poles, this should result in an order of magnitude increase in dampening for the closed-loop system.
	
	For the choice of $C(s)$ in (\ref{ch example: PID controller C(s)}), our augmented-system has the state space realisation
	\begin{align}
		\dot{x} &= Ax + B_1w + B_2u, \nonumber\\
		z       &= C_1x, \label{ch example: math: PID augmented system realisation}
	\end{align}
	where
	\begin{align*}
		A &= \begin{bmatrix}
			\begin{tabular}{ c c c} 	
				$0$ & $1$ & $0$ \\ 
				$0$ & $0$ & $1$ \\
				$3203$ & $-80511.58$ & $-12.812$
			\end{tabular}
		\end{bmatrix}, \\ B_2 &= \begin{bmatrix}
			\begin{tabular}{ c } 	
				$0$ \\
				$0$ \\
				$1$
			\end{tabular}
		\end{bmatrix}, \quad B_1 = \begin{bmatrix}
			\begin{tabular}{ c } 	
				$0$ \\
				$0$ \\
				$-1$
			\end{tabular}
		\end{bmatrix}, \\\\
		C_1 &= \begin{bmatrix}
			\begin{tabular}{ c c c} 	
				$-3203$ & $-448.42$ & $-12.812$ 
			\end{tabular}
		\end{bmatrix}.
	\end{align*}
	
	In order to apply SNI state feedback to the augmented system (\ref{ch example: math: PID augmented system realisation}), we form the matrices $Q$, $A_q$ from (\ref{math: Schur Q definition}) and (\ref{math: A_q equation}) respectively.
	
	\begin{align}
		Q &= I-B_2(C_1B_2)^{-1}C_1 = \begin{bmatrix}
			\begin{tabular}{ c c c} 	
				$1$ & $0$ & $0$ \\ 
				$0$ & $1$ & $0$ \\
				$-250$ & $-35$ & $0$
			\end{tabular}
		\end{bmatrix}, \nonumber\\
		A_q &= QA = \begin{bmatrix}
			\begin{tabular}{ c c c} 	
				$0$ & $1$ & $0$ \\ 
				$0$ & $0$ & $1$ \\
				$0$ & $-250$ & $-35$
			\end{tabular}
		\end{bmatrix}.
	\end{align}
	The eigenvalues of $A_q$ are $\sigma(A_q) = \{0,-10,-25\}$.
	We are now free to use Theorem~\ref{theorem: SNI ARE synthesis theorem no A22} in order to synthesise a controller that results in a closed-loop SNI system with a prescribed degree of stability. From Theorem~\ref{chap how far pert: distinct eigenvalues no A22}, we note that we can choose a value of $\epsilon$ within the range $0<\epsilon<10$ before we need to consider a dimension change in the solution of the ARE used in the controller synthesis. For the sake of this example a prescribed degree of stability of $\epsilon=9.5$ is chosen.
	
	Theorem~\ref{theorem: SNI ARE synthesis theorem no A22} states that the state feedback controller that will result in an SNI closed-loop system with a prescribed degree of stability $\epsilon$, is given by
	\begin{align} 
		K &= -(C_1B_2)^{-1}(C_1A + \epsilon C_1 ) \nonumber \\
		&= \begin{bmatrix}
			\begin{tabular}{ c c c} 	
				$-5578$ & $79929.08$ & $-57.312$ 
			\end{tabular}
		\end{bmatrix}. \label{ch examp: math: k matrix reduced}
	\end{align} 
	Applying (\ref{ch examp: math: k matrix reduced}) to the system (\ref{ch example: math: PID augmented system realisation}) using the state feedback control law $u=Kx$, the corresponding closed-loop uncertain system has a state space representation
	\begin{align}
		\dot{x} &= (A + B_2K)x + B_1w, \nonumber\\
		z       &= C_1x ,       \label{ch example: math: closed-loop system}
	\end{align}
	where
	\begin{align}
		A_{cl} &= A + B_2K = \begin{bmatrix}
			\begin{tabular}{ c c c} 	
				$0$ & $1$ & $0$ \\ 
				$0$ & $0$ & $1$ \\
				$-2375$ & $-582.5$ & $-44.5$
			\end{tabular}
		\end{bmatrix}, \\
		B_2 &= \begin{bmatrix}
			\begin{tabular}{ c } 	
				$0$ \\
				$0$ \\
				$1$
			\end{tabular}
		\end{bmatrix}, \quad B_1 = \begin{bmatrix}
			\begin{tabular}{ c } 	
				$0$ \\
				$0$ \\
				$-1$
			\end{tabular}
		\end{bmatrix}, \nonumber\\\\
		C_1 &= \begin{bmatrix}
			\begin{tabular}{ c c c} 	
				$-3203$ & $-448.42$ & $-12.812$ 
			\end{tabular}
		\end{bmatrix}.
	\end{align}
	The closed-loop poles are located at $\sigma(A_{cl}) =\{-25,-10,$ $-9.5\}$ and the closed-loop system is stable.
	
\begin{figure}[!ht]
		\centering
		\includegraphics[width=3in]{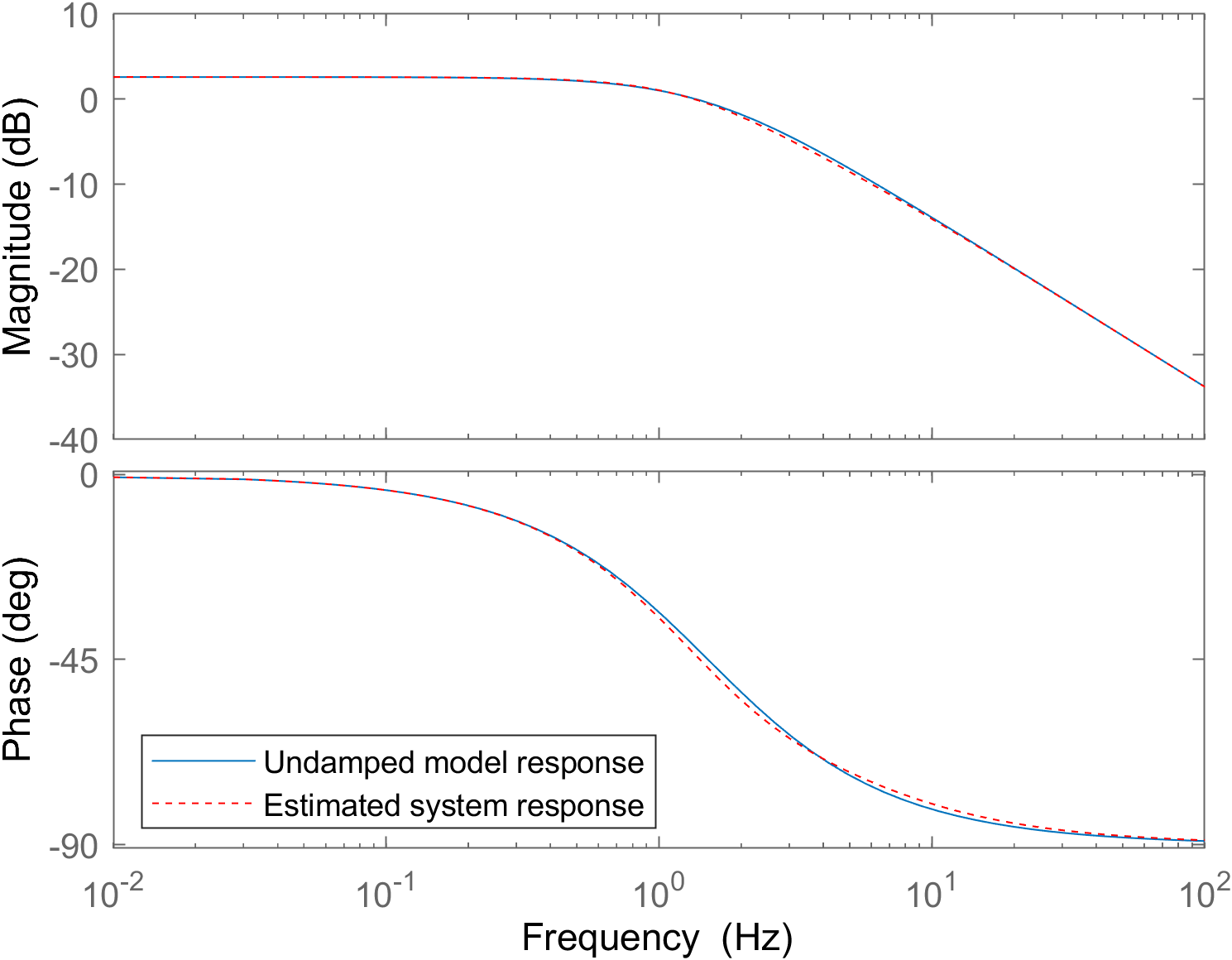}
		\caption{The frequency response of the SNI closed-loop transfer function (\ref{ch example: PID example: closed-loop tf}) and approximated model (\ref{ch exam: math: lightly damped GN}).}
		\label{ch examp: fig: closed-loop freq response of tf}
	\end{figure}

	The corresponding closed-loop transfer function is 
	\begin{align}
		G_{cl}(s) &= C_1(sI - A - B_2K)^{-1}B_1 \nonumber\\
		&= \frac{12.81 s^2 + 448.4 s + 3203}{s^3 + 44.5 s^2 + 582.5 s + 2375}, \label{ch example: PID example: closed-loop tf}
	\end{align}
	and has the frequency response shown in Figure~\ref{ch examp: fig: closed-loop freq response of tf}. Figure~\ref{ch examp: fig: closed-loop freq response of tf} shows the frequency response of the closed-loop system for both the approximated model (\ref{ch exam: math: lightly damped GN}) and the actual system response of the simulated beam. Observe that both responses are SNI. The uncertainty is lossless NI and we can easily verify that the gain conditions of Lemma~\ref{lanzon 2017 stability result} are satisfied. Thus, the overall system is internally stable.
	
	\section{Conclusion}
	
	In this paper we have presented a method of augmenting an NI system with an additional controller so that the open loop system satisfies the necessary assumptions to apply SNI state feedback. We demonstrated two methods for modelling NI systems with a robust control framework using additive and multiplicative uncertainty and a sum of second order transfer functions corresponding to the resonant modes of the system. It was shown that both modelling methods have orthogonal input and output matrices when considering their corresponding state space realisations. As a result, neither method satisfies the assumptions needed to apply many of the NI feedback methods found in the literature. By augmenting these models with an additional integral or PID controller, the open loop dynamics were sufficiently altered in order to overcome this limitation. Moreover, the affects of the two controllers on the degree of stability that the closed-loop system could achieve under SNI state feedback was quantified. Finally, an illustrative example demonstrated how our approach could be applied to a simulated flexible cantilever. In Futures works the results presented here could be extended to the output feedback case and extended to time-varying systems.   
	
	\bibliography{references}             

\begin{thebibliography}{4}
\providecommand{\natexlab}[1]{#1}
\providecommand{\url}[1]{\texttt{#1}}
\providecommand{\urlprefix}{URL }
\expandafter\ifx\csname urlstyle\endcsname\relax
  \providecommand{\doi}[1]{doi:\discretionary{}{}{}#1}\else
  \providecommand{\doi}{doi:\discretionary{}{}{}\begingroup
  \urlstyle{rm}\Url}\fi

\bibitem[{Able(1956)}]{Abl:56}
Able, B. (1956).
\newblock Nucleic acid content of microscope.
\newblock \emph{Nature}, 135, 7--9.

\bibitem[{Able et~al.(1954)Able, Tagg, and Rush}]{AbTaRu:54}
Able, B., Tagg, R., and Rush, M. (1954).
\newblock Enzyme-catalyzed cellular transanimations.
\newblock In A.~Round (ed.), \emph{Advances in Enzymology}, volume~2, 125--247.
  Academic Press, New York, 3rd edition.

\bibitem[{Keohane(1958)}]{Keo:58}
Keohane, R. (1958).
\newblock \emph{Power and Interdependence: World Politics in Transitions}.
\newblock Little, Brown \& Co., Boston.

\bibitem[{Powers(1985)}]{Pow:85}
Powers, T. (1985).
\newblock Is there a way out?
\newblock \emph{Harpers}, 35--47.

\end{thebibliography}


\begin{thebibliography}{16}
\providecommand{\natexlab}[1]{#1}
\providecommand{\url}[1]{\texttt{#1}}
\providecommand{\urlprefix}{URL }
\expandafter\ifx\csname urlstyle\endcsname\relax
  \providecommand{\doi}[1]{doi:\discretionary{}{}{}#1}\else
  \providecommand{\doi}{doi:\discretionary{}{}{}\begingroup
  \urlstyle{rm}\Url}\fi

\bibitem[{Bernstein(2009)}]{Bernstein2009}
Bernstein, D.S. (2009).
\newblock \emph{Matrix mathematics : theory, facts, and formulas}.
\newblock Princeton University Press, Princeton, New Jersey, USA.

\bibitem[{Bhowmick and Patra(2020)}]{10.1016/j.automatica.2019.108735}
Bhowmick, P. and Patra, S. (2020).
\newblock Solution to negative-imaginary control problem for uncertain lti
  systems with multi-objective performance.
\newblock \emph{Automatica}, 112(C).
\newblock \doi{10.1016/j.automatica.2019.108735}.

\bibitem[{Dannatt et~al.(2020)Dannatt, Petersen, and
  Lanzon}]{dannatt2020strictly}
Dannatt, J., Petersen, I.R., and Lanzon, A. (2020).
\newblock Strictly negative imaginary state feedback control with a prescribed
  degree of stability.
\newblock \emph{Automatica}, 119, 109079.

\bibitem[{Hu et~al.(2022)Hu, Lennox, and
  Arvin}]{10.1016/j.automatica.2022.110235}
Hu, J., Lennox, B., and Arvin, F. (2022).
\newblock Robust formation control for networked robotic systems using negative
  imaginary dynamics.
\newblock \emph{Automatica}, 140(C).
\newblock \doi{10.1016/j.automatica.2022.110235}.
\newblock \urlprefix\url{https://doi.org/10.1016/j.automatica.2022.110235}.

\bibitem[{Lanzon and Chen(2017)}]{lanzon2017feedback}
Lanzon, A. and Chen, H.J. (2017).
\newblock Feedback stability of negative imaginary systems.
\newblock \emph{IEEE Transactions on Automatic Control}, 62(11), 5620--5633.

\bibitem[{Lanzon and Petersen(2008)}]{Lanzon2008}
Lanzon, A. and Petersen, I.R. (2008).
\newblock Stability robustness of a feedback interconnection of systems with
  negative imaginary frequency response.
\newblock \emph{IEEE Transactions on Automatic Control}, 53(4), 1042--1046.
\newblock \doi{10.1109/TAC.2008.919567}.

\bibitem[{Liu and Xiong(2016)}]{LIU201647}
Liu, M. and Xiong, J. (2016).
\newblock On non-proper negative imaginary systems.
\newblock \emph{Systems and Control Letters}, 88, 47--53.
\newblock \doi{https://doi.org/10.1016/j.sysconle.2015.10.008}.

\bibitem[{Mabrok et~al.(2015)Mabrok, Kallapur, Petersen, and
  Lanzon}]{Mabrok2015}
Mabrok, M., Kallapur, A.G., Petersen, I.R., and Lanzon, A. (2015).
\newblock A generalized negative imaginary lemma and {Riccati}-based static
  state-feedback negative imaginary synthesis.
\newblock \emph{Systems {\&} Control Letters}, 77, 63--68.
\newblock \doi{10.1016/j.sysconle.2015.01.008}.

\bibitem[{Mabrok et~al.(2012)Mabrok, Kallapur, Petersen, and
  Lanzon}]{Mabrok2012a}
Mabrok, M.A., Kallapur, A.G., Petersen, I.R., and Lanzon, A. (2012).
\newblock {Stabilization of uncertain negative-imaginary systems using a
  Riccati equation approach}.
\newblock In \emph{Proceedings of the 2012 First International Conference on
  Innovative Engineering Systems}, 255--259. IEEE, Alexandria, Egypt.
\newblock \doi{10.1109/ICIES.2012.6530879}.

\bibitem[{Petersen(2015)}]{Petersen2015}
Petersen, I.R. (2015).
\newblock {Physical interpretations of negative imaginary systems theory}.
\newblock In \emph{2015 10th Asian Control Conference (ASCC)}, 1--6. IEEE, Kota
  Kinabalu, Malaysia.
\newblock \doi{10.1109/ASCC.2015.7244460}.

\bibitem[{Petersen(2016)}]{petersen2016negative}
Petersen, I.R. (2016).
\newblock Negative imaginary systems theory and applications.
\newblock \emph{Annual Reviews in Control}, 42, 309--318.

\bibitem[{Petersen and Lanzon(2010)}]{Petersen2010}
Petersen, I.R. and Lanzon, A. (2010).
\newblock Feedback control of negative-imaginary systems.
\newblock \emph{IEEE Control Systems Magazine}, 30(5), 54--72.
\newblock \doi{10.1109/MCS.2010.937676}.

\bibitem[{Preumont(1997)}]{preumont1997vibration}
Preumont, A. (1997).
\newblock \emph{Vibration control of active structures}, volume~2.
\newblock Springer.

\bibitem[{Ren et~al.(2021)Ren, Xiong, and Ho}]{REN2021109157}
Ren, D., Xiong, J., and Ho, D.W. (2021).
\newblock Static output feedback negative imaginary controller synthesis with
  an {$H_\infty$} norm bound.
\newblock \emph{Automatica}, 126, 109157.
\newblock \doi{https://doi.org/10.1016/j.automatica.2020.109157}.

\bibitem[{Reyes(2019)}]{reyes2019controller}
Reyes, G.P.S. (2019).
\newblock \emph{Controller Synthesis for Strictly Negative Imaginary Systems
  Via Riccati Equations and Linear Matrix Inequalities}.
\newblock The University of Manchester (United Kingdom).

\bibitem[{Song et~al.(2012)Song, Lanzon, Patra, and Petersen}]{Song2012}
Song, Z., Lanzon, A., Patra, S., and Petersen, I.R. (2012).
\newblock {A negative-imaginary lemma without minimality assumptions and robust
  state-feedback synthesis for uncertain negative-imaginary systems}.
\newblock \emph{Systems {\&} Control Letters}, 61(12), 1269--1276.
\newblock \doi{10.1016/j.sysconle.2012.08.002}.

\end{thebibliography}
	
\end{document}